\documentclass[a4paper]{article}
\usepackage{amssymb}
\usepackage{amsfonts}
\usepackage{amsmath}
\usepackage{amsthm}
\usepackage{authblk}
\usepackage[USenglish]{babel}
\usepackage{esint}
\usepackage{graphicx}
\usepackage[colorlinks=true]{hyperref}
\usepackage[latin1]{inputenc}
\newtheorem{theorem}{Theorem}[section]
\newtheorem{lemma}[theorem]{Lemma}
\newtheorem{proposition}[theorem]{Proposition}
\newtheorem{corollary}[theorem]{Corollary}
\newtheorem{remark}[theorem]{Remark}
\newtheorem{definition}[theorem]{Definition}

\numberwithin{equation}{section}

\newcommand{\R}{\mathbb R}

\newcommand{\Di}{\mathbb D}

\newcommand{\I}{\mathbb I}

\newcommand{\mrm}{\mathrm}

\renewcommand{\a}{\alpha}

\newcommand{\D}{\Delta}

\newcommand{\la}{\lambda}

\newcommand{\s}{\sigma}

\renewcommand{\O}{\Omega}
\newcommand{\wt}{\widetilde}

\newcommand{\ol}{\overline}

\newcommand{\ub}{\underbrace}
\newcommand{\fr}{\frac}
\newcommand{\pa}{\partial}
\newcommand{\n}{\nabla}

\newcommand{\sm}{\setminus}

\newcommand{\sub}{\subset}

\newcommand{\eq}{\equiv}
\newcommand{\ox}{\otimes}

\renewcommand{\c}{\circ}
\newcommand{\cd}{\cdot}
\newcommand{\ds}{\dots}

\newcommand{\tx}{\text}
\newcommand{\q}{\quad}
\renewcommand{\l}{\left}
\renewcommand{\r}{\right}

\newcommand{\bthm}{\begin{theorem}}
\newcommand{\ethm}{\end{theorem}}
\newcommand{\blem}{\begin{lemma}}
\newcommand{\elem}{\end{lemma}}
\newcommand{\bprop}{\begin{proposition}}
\newcommand{\eprop}{\end{proposition}}
\newcommand{\bcor}{\begin{corollary}}
\newcommand{\ecor}{\end{corollary}}
\newcommand{\brem}{\begin{remark}}
\newcommand{\erem}{\end{remark}}
\newcommand{\bdefi}{\begin{definition}}
\newcommand{\edefi}{\end{definition}}
\newcommand{\bpf}{\begin{proof}}
\newcommand{\epf}{\end{proof}}
\newcommand{\bl}{\begin{array}{l}}
\newcommand{\bll}{\begin{array}{ll}}
\newcommand{\barr}{\begin{array}}
\newcommand{\earr}{\end{array}}
\newcommand{\bite}{\begin{itemize}}
\newcommand{\eite}{\end{itemize}}
\newcommand{\bequ}{\begin{equation}}
\newcommand{\eequ}{\end{equation}}
\newcommand{\beqa}{\begin{eqnarray}}
\newcommand{\eeqa}{\end{eqnarray}}
\newcommand{\beqy}{\begin{eqnarray*}}
\newcommand{\eeqy}{\end{eqnarray*}}

\allowdisplaybreaks

\begin{document}

\everymath{\displaystyle}

\title{Uniform bounds for solutions to elliptic problems on simply connected planar domains}
\author{Luca Battaglia\thanks{Universit\`a degli Studi Roma Tre, Dipartimento di Matematica e Fisica, Largo S. Leonardo Murialdo 1, 00146 Roma - lbattaglia@mat.uniroma3.it}}
\date{}

\maketitle\

\begin{abstract}
\noindent
We consider the following elliptic problems on simply connected planar domains
$$\l\{\bll-\D u=\la|x|^{2\a}K(x)e^u&\tx{in }\O\\u=0&\tx{on }\pa\O\earr\r.;\q\q\q\l\{\bll-\D u=|x|^{2\a}K(x)u^p&\tx{in }\O\\u>0&\tx{in }\O\\u=0&\tx{on }\pa\O\earr\r.;$$
with $\a>-1,\la>0,p>1,0<K(x)\in C^1\l(\ol\O\r)$.\\
We show that any solution to each problem must satisfy a uniform bound on the mass, which is given respectively by $\la\int_\O|x|^{2\a}K(x)e^u\mrm dx$ and $p\int_\O|x|^{2\a}K(x)u^{p+1}\mrm dx$. The same results applies to some systems and more general non-linearities.\\
The proofs are based on the Riemann mapping theorem and a Poho\v zaev-type identity.
\end{abstract}\

\section{Introduction}\

We are interested in the following PDE, known as \emph{Liouville equation}:
\bequ\label{lioueq}
\l\{\bll-\D u=\la|x|^{2\a}K(x)e^u&\tx{in }\O\\u=0&\tx{on }\pa\O\earr\r.,
\eequ
with $\O\ni0$ being a smooth bounded planar domain, $\a>-1,\la>0$ and $0<K(x)\in C^1\l(\O\r)$.\\
Equation \eqref{lioueq} has been very deeply studied in the last decades due to its applications in geometry and physics. It may be considered as a \emph{critical} elliptic problem on planar domain, as the exponential nonlinearity is a natural counterpart of the Sobolev critical exponent in dimension greater or equal than $3$.\\
Solutions to \eqref{lioueq} can been found either variationally (\cite{dja,djamal,bdm,cm}) or by computing the Leray-Schauder degree (\cite{cl03,cl13}), and blowing-up families have also been constructed (\cite{egp,dkm}). In all of these cases the geometry and topology of the domain $\O$ play a fundamental role. See also \cite{bp} for existence of blowing-up solutions to \eqref{lioueq}.\\
In this paper we give a \emph{mass bound} for solutions to \eqref{lioueq} when $\O$ is simply connected, namely we show that any solution must satisfy a uniform bound on the $L^1$ norm of the laplacian $\rho:=\la\int_\O|x|^{2\a}K(x)e^u$. Such a quantity plays an important role especially in the variational formulation of the problem and sometimes, to stress its importance, it is used as a parameter in place of $\la$, with the equation in \eqref{lioueq} rewritten as $-\D u=\rho\fr{|x|^{2\a}K(x)e^u}{\int_\O|x|^{2\a}K(x)e^u\mrm dx}$.\\
The following results extends a previous one on the unit disk (\cite{barmal}, Proposition $5.7$).\\

\bthm\label{liou}${}$\\
Let $\O$ be a simply connected planar domain and $u$ be a solution to \eqref{lioueq}.\\
Then, there exists $\rho_0>0$, not depending on $u$, such that $\la\int_\O|x|^{2\a}K(x)e^u\mrm dx\le\rho_0$.
\ethm\

We are also considering the \emph{Hénon-Lane-Emden} equation:
\bequ\label{hleq}
\l\{\bll-\D u=|x|^{2\a}K(x)u^p&\tx{in }\O\\u>0&\tx{in }\O\\u=0&\tx{on }\pa\O\earr\r.,
\eequ
with $p>1$ and $\O,\a,K(x)$ as before. The power-type nonlinearity in \eqref{hleq} is \emph{subcritical} on planar domains for any $p$, hence positive solutions can be easily found. Nonetheless, it is interesting to investigate the asymptotic behavior of solutions as the exponent $p$ goes to $+\infty$.\\
Despite the different structure, the latter problem shares surprising similarities with blow-up analysis for equation \eqref{lioueq} (see \cite{dip,dgip}), and in particular the problem heavily depends on the shape of $\O$. See also \cite{rw1,rw2} for asymptotic analysis of solutions to \eqref{hleq}.
In the regular case $\a=0$ solutions to \eqref{hleq} have been found on multiply connected domains with arbitrarily large values of the mass, which in this case is given by $p\int_\O|x|^{2\a}K(x)u^{p+1}\mrm dx$ (see \cite{emp}). On the other hand, when $\O$ is convex (\cite{gt,dgip}) or strictly star-shaped (\cite{ks}) different bounds on the mass have been given, which are equivalent to giving an upper bound on the number of blow-up points.\\
Here we fill the gap between the two results by showing that, in the same spirit as Theorem \ref{liou}, similar bounds hold true for solutions to \eqref{hleq} on simply connected domains.\\

\bthm\label{hle}${}$\\
Let $\O$ be a simply connected planar domain, $p_0>1$ and $u$ be a solution to \eqref{hleq} with $p\ge p_0$.\\
Then, there exists $\rho_0>0$, not depending on $u$ nor on $p$, such that $p\int_\O|x|^{2\a}K(x)u^{p+1}\mrm dx\le\rho_0$.
\ethm\

We will also provide results similar to Theorems \ref{liou}, \ref{hle} to some \emph{Liouville systems}, namely systems of PDEs with the same features as \eqref{lioueq}. Such problems have been increasingly studied in the last years, especially in the case when the matrix of coefficients $A$ is a Cartan matrix of some Lie algebra.\\
To get a mass bound for solutions on simply connected domains, we need the matrix $A$ to be positive definite, which in the case of Cartan matrices holds true. Such a result had already been proven when $\O$ is the unit disk, in \cite{batmal} (Theorem $1.3$) for the $SU(3)$ Toda system and in the author's PhD thesis \cite{bat0} for general systems. A similar estimate was proved also in \cite{bar} for general systems on strictly star-shaped domains.\\

\bthm\label{lious}${}$\\
Let $\O$ be a simply connected planar domain, $A=\{a_{ij}\}_{i,j=1}^N$ be a positive definite matrix and $u=(u_1,\ds,u_N)$ be a solution to
\bequ\label{liouseq}
\l\{\bll-\D u_i=\sum_{j=1}^Na_{ij}\la_j|x|^{2\a_j}K_j(x)e^{u_j}&\tx{in }\O\\u_i=0&\tx{on }\pa\O\earr\r.\q\q\q i=1,\ds,N,
\eequ
with $\O\ni0,\a_i>-1,\la_i>0,0<K_i(x)\in C^1\l(\ol\O\r)$.\\
Then, there exists $\rho_0>0$, not depending on $u$, such that $\la_i\int_\O|x|^{2\a_i}K_i(x)e^{u_i}\mrm dx\le\rho_0$ for all $i$'s.
\ethm\
 
Finally, similar estimates also hold true for some more general critical linearities.\\
Roughly speaking, we need a positive potential $W(x)$ to be not too singular at $0$ and the non-linearity $F'(u)$ to grow at least as fast as its anti-derivative $F(u)$. This case includes most exponential functions, including $F'(u)=\la u^{p-1}e^{u^p}$ with $1<p<2$, which was studied in \cite{dm1,dgm,dm2}.\\

\bthm\label{gen}${}$\\
Let $\O$ be a simply connected planar domain and $u$ be a solution to
\bequ\label{geneq}
\l\{\bll-\D u=W(x)F'(u)&\tx{in }\O\\u=0&\tx{on }\pa\O\earr\r.,
\eequ
with $W(x),F(u)$ satisfying
\bequ\label{cond}
\l\{\bll W(x)\in C^1_{\mrm{loc}}\l(\ol\O\sm\{0\}\r)\\0\le W(x)\le C|x|^{2\a}&\a>-1\\|\n W(x)||x|\le CW(x)\earr\r.\q\q\q\q\q\q\l\{\bl F\in C^1(\R)\\0\le F'(u)\le Ce^{Cu^2}\\F(u)\le C(1+F'(u))\earr\r..
\eequ
Then, there exists $\rho_0>0$, not depending on $u$, such that $\int_\O W(x)F'(u)\mrm dx\le\rho_0$.
\ethm\

It is interesting to compare all these results with previous works concerning existence and non-existence of solutions. This is the content of the following remark.\\

\brem${}$\\
\bite
\item If the domain $\O$ is not simply connected, then each of the problem we are considering can have solutions with arbitrarily high values of the mass. This was done in \cite{dkm,dja,djamal,bdm,cm,car}  for problem \eqref{lioueq}, in \cite{emp,epw} for \eqref{hleq}, in \cite{bjmr,bat1,bat2} for some systems of the type \eqref{liouseq} and in \cite{dm1,dgm,dm2} for some nonlinearities of the type \eqref{geneq}.\\
\item In \cite{egp} the authors prove that for any $M>0$ there exist a simply connected \emph{dumbbell}-shaped domain $\O_M$ and a solution to \eqref{lioueq} on $\O_M$ with $\la\int_{\O_M}|x|^{2\a}K(x)e^u\mrm dx\ge M$; the same argument also works for problems \eqref{hleq}, \eqref{geneq} (see \cite{epw} and \cite{dm1}, respectively). The results presented here complement the latter, since Theorem \ref{liou} implies that there cannot exists any $\O_M$ such that the property shown in \cite{egp} holds for any $M$.\\
\item In Theorem \ref{lious} it is essential to assume the matrix $A$ to be positive definite. Otherwise, in \cite{gp,mpw,pr,bp} the authors build solutions to \eqref{liouseq} whose masses can be arbitrarily large also on simply connected domains.\\
\item If one allows more than one singularity, namely replaces the singular term $|x|^{2\a}$ with $|x|^{2\a}\prod_{i=1}^N|x-x_i|^{2\a_i}$ for some $x_i\in\O\sm\{0\},\a_i>-1$, then uniform mass bounds do not seem to be true anymore. In fact, in this case Theorem $1.4$ in \cite{cl13} shows that the Leray-Schauder degree of \eqref{lioueq} does not vanish for arbitrarily high values of the mass.\\
\item Finally, assuming $K(x)$ to be positive is essential. In fact, in \cite{lr,dl,dlr} the authors show existence of solutions to the Liouville equation \eqref{lioueq} with sign-changing potential even in the case of simply connected domain; here, a crucial role seems to be played not by the topology of $\O$ but rather of the set $\{x\in\O:\,K(x)>0\}$.\\
\eite
\erem\

The main tools to prove Theorems \ref{liou}, \ref{hle}, \ref{lious}, \ref{gen} will be the Riemann Mapping Theorem and a Poho\v zaev-type identity.\\
We will recall these very well-known results in Section $2$, as well as some other preliminary. Then, in Section $3$ we will prove the main results of this paper.\\

\section{Preliminaries}\

Let us recall some facts which will be used in the proof of the results of this paper.\\
We start with a very classical and powerful tool, the Riemann Mapping Theorem. Such a results will allow to conformally deform the simply connected domain $\O$ into the unit disk; in such a way, the PDE defined on $\O$ is transformed into a new equation on the disk, different from the original but with similar features.\\
We actually need a refined version of the theorem by Carathéodory, which ensures that the conformal factor appearing in the new PDE is not singular on the boundary of $\O$.\\

\bthm[Riemann Mapping Theorem, Carathéodory's Theorem]\label{map}${}$\\
Let $\O\ni0$ be a smooth simply connected planar domain and $\Di\sub\R^2$ be the unit disk.\\
Then, there exists a conformal diffeomorphism $\Phi:\O\to\Di$, smooth up to $\pa\O$, such that $\Phi(0)=0$.\\
Moreover, if $u$ solves $\l\{\bll-\D u=f(x,u)&\tx{in }\O\\u=0&\tx{on }\pa\O\earr\r.$, then $v:=u\c\Phi^{-1}$ solves\\$\l\{\bll-\D v=\fr{1}{\det\l(D\Phi(y)\r)}f\l(\Phi^{-1}(y),v\r)&\tx{in }\Di\\v=0&\tx{on }\pa\Di\earr\r.$.
\ethm\

Our proofs will also use the Poho\v zaev identity, an often-used instrument to show non-existence of solutions to elliptic PDEs.\\
Such a result is usually stated for solutions having at least a $W^{2,2}$ regularity, which in general does not hold true for solutions to \eqref{lioueq}, \eqref{hleq}, \eqref{liouseq}, \eqref{geneq} if $\a$ is negative. Anyway, in the proof of the theorems we will verify that we are still in position to apply the following result.\\

\bthm[Poho\v zaev Identity]\label{poho}${}$\\
Let $u$ be a sufficiently regular solution to
$$\l\{\bll-\D v=\pa_vG(y,v)&\tx{in }\Di\\v=0&\tx{on }\pa\Di\earr\r..$$
Then, it satisfies
$$\fr{1}2\int_{\pa\Di}(\n v\cd\nu(y))^2\mrm d\s(y)=2\int_\Di G(y,v)\mrm dy+\int_\Di(\n_yG(y,v)\cd y)\mrm dy-\int_{\pa\Di}G(y,v)\mrm d\s(y).$$
If $v=(v_1,\ds,v_N)$ solves 
$$\l\{\bll-\D v_i=\sum_{j=1}^Na_{ij}\pa_{v_j}G_j(y,v_j)&\tx{in }\Di\\v_i=0&\tx{on }\pa\Di\earr\r.\q\q\q i=1,\ds,N,$$
with $A=\{a_{ij}\}_{i,j=1,\ds,N}$ being a non-singular matrix, then $u$ satisfies
$$\fr{1}2\sum_{i,j=1}^Na^{ij}\int_{\pa\Di}(\n v_i\cd\nu(y))(\n v_j\cd\nu(y))\mrm d\s(y)=\sum_{i=1}^N\l(2\int_\Di G_i(y,v_i)\mrm dy+\int_\Di(\n_yG_i(y,v_i)\cd y)\mrm dy-\int_{\pa\Di}G_i(y,v_i)\mrm d\s(y)\r).$$
\ethm\

We finally need some \emph{a priori} estimates for solutions to \eqref{hleq}, which are essential to adapt the argument for \eqref{lioueq}.\\
The following result was originally stated in \cite{ks} for the case $K\eq1,\a=0$ but the same argument, based on estimates from \cite{dln} and the celebrated moving plane technique from \cite{gnn}, seems to be working more generally.\\

\bthm\label{ks}${}$\\
Let $u$ be a solution to \eqref{hleq} with $p\ge p_0$.\\
Then, there exists $C_0>0$, not depending on $p$ nor $u$, such that $\sup_\O u=\|u\|_{L^\infty(\O)}\le C_0$.
\ethm\

\section{Proofs}\

We are now in position to prove the results stated in the introduction.\\
Since all the proofs are rather similar to each other, we will give more details for Theorem \ref{liou} but we will skip some for the other theorems.\\

\bpf[Proof of Theorem \ref{liou}]${}$\\
Let $u$ be a solution to \eqref{lioueq} on $\O$ and $\Phi:\O\to\Di$ be the Riemann mapping described in Theorem \ref{map}. Then, $v:=u\c\Phi^{-1}$ will solve $\l\{\bll-\D v=\la\wt K(y)|y|^{2\a}e^v&\tx{in }\Di\\u=0&\tx{on }\pa\Di\earr\r.$ for some $0<\wt K(y)\in C^1\l(\ol\Di\r)$.\\
We want to apply to $v$ Theorem \ref{poho} with $G(y,v)=\la|y|^{2\a}\wt K(y)e^v$. If $\a>0$ this is immediate since standard regularity gives $v\in C^{2,2\a}\l(\ol\Di\r)$, but in case $\a<0$ we only have $v\in W^{2,q}(\Di)$ with $q<\fr{1}{-\a}$, therefore we need an \emph{ad hoc} argument. Poho\v zaev identity is based on applying the divergence theorem to $(\n v\cd y)\n v-\fr{|\n v|^2}2y$ and $\la\wt K(y)|y|^{2\a}e^uy$, so we need to check that both vector fields are in $W^{1,1}\l(\Di\r)$.\\
Concerning the former field, we have
$$D\l((\n v\cd y)\n v-\fr{|\n v|^2}2y\r)=(\n v\cd y)D^2v+\n v\ox\n v+\l(D^2v,y\r)\ox\n v-\fr{|\n v|^2}2\I_2-\l(D^2v,\n v\r)\ox y;$$
since we already know that $|\n v|^2\in L^1(\Di)$, we suffice to check that $\l|D^2 v\r||\n v||y|\in L^1(\Di)$. We have $D^2v\in L^1(\Di)$ and moreover, since $|\D v|\le C|y|^{2\a}$, by the Green's representation formula we deduce $|\n v|\le C\int_\Di\fr{|\eta|^{2\a}}{|y-\eta|}\mrm d\eta\le\fr{C}{|y|}$, hence $|\n v||y|\in L^\infty(\Di)$ and we are done.\\
The other vector field verifies
$$D\l(\la\wt K(y)|y|^{2\a}e^uy\r)=\la K(y)|y|^{2\a}e^v\I_2+\la|y|^{2\a}e^v\l(\n\wt K(y)\ox y\r)+2\a\la|y|^{2\a-2}\wt K(y)e^v(y\ox y),$$
which is in $L^1(\Di)$ because each term can be estimated by constant times $|y|^{2\a}$.\\
We are therefore in position to use Theorem \ref{poho}, which gives:
$$\fr{1}2\int_{\pa\Di}(\n v\cd\nu(y))^2\mrm d\s(y)=2\la\int_\Di|y|^{2\a}\wt K(y)e^v\mrm dy+\la\int_\Di\l(2\a|y|^{2\a}\wt K(y)+|y|^{2\a}\l(\n\wt K(y)\cd y\r)\r)e^v\mrm dy-\la\int_{\pa\Di}\wt K(y)e^v\mrm d\s(y).$$
On the left-hand side we can use H\"older's inequality and integrate by parts:
$$\fr{1}2\int_{\pa\Di}(\n v\cd\nu(y))^2\mrm d\s(y)\ge\fr{1}{4\pi}\l(\int_{\pa\Di}\n v\cd\nu(y)\mrm d\s(y)\r)^2=\fr{1}{4\pi}\l(\int_\Di\D v\mrm dy\r)^2=\fr{1}{4\pi}\l(\la\int_\Di|y|^{2\a}\wt K(y)e^v\mrm dy\r)^2;$$
on the right-hand side we exploit the positivity of $\wt K(y)$ and the boundedness of $\n\wt K(y)$:
\beqy
&&2\la\int_\Di|y|^{2\a}\wt K(y)e^v\mrm dy+\la\int_\Di\l(2\a|y|^{2\a}\wt K(y)+|y|^{2\a}\l(\n\wt K(y)\cd y\r)\r)e^v\mrm dy-\la\int_{\pa\Di}\wt K(y)e^v\mrm d\s(y)\\
&\le&2(1+\a)\la\int_\Di|y|^{2\a}\wt K(y)e^v\mrm dy+\la\int_\Di|y|^{2\a}\l(\n\wt K(y)\cd y\r)e^v\mrm dy\\
&\le&\ub{\l(2(1+\a)+\fr{\sup_{y\in\Di}\l|\n\wt K(y)\r|}{\inf_{y\in\Di}\wt K(y)}\r)}_{=:\fr{\rho_0}{4\pi}}\l(\la\int_\Di|y|^{2\a}\wt K(y)e^v\mrm dy\r).
\eeqy
Putting the two estimates together we get
$$\fr{1}{4\pi}\l(\la\int_\Di|y|^{2\a}\wt K(y)e^v\mrm dy\r)^2\le\fr{\rho_0}{4\pi}\l(\la\int_\Di|y|^{2\a}\wt K(y)e^v\mrm dy\r),$$
hence we conclude
$$\la\int_\O|x|^{2\a}K(x)e^u\mrm dx=\la\int_\Di|y|^{2\a}\wt K(y)e^v\mrm dy\le\rho_0.$$
\epf\

\bpf[Proof of Theorem \ref{hle}]${}$\\
As in the proof of Theorem \ref{liou}, we take a solution to \eqref{hleq} and transform it, via the Riemann mapping, into a solution to $\l\{\bll-\D v=|y|^{2\a}\wt K(y)v^p&\tx{in }\Di\\v>0&\tx{in }\Di\\v=0&\tx{on }\pa\Di\earr\r.,$ with some $0<\wt K(y)\in C^1\l(\ol\Di\r)$.\\
By the same argument as in Theorem \ref{liou}, we are allowed to apply Theorem \ref{poho}, this time with $G(y,v)=|y|^{2\a}\wt K(y)\fr{v^{p+1}}{p+1}$, which reads as
$$\fr{1}2\int_{\pa\Di}(\n v\cd\nu(y))^2\mrm d\s(y)=2\int_\Di|y|^{2\a}\wt K(y)\fr{v^{p+1}}{p+1}\mrm dy+\int_\Di\l(2\a|y|^{2\a}\wt K(y)+|y|^{2\a}\l(\n\wt K(y)\cd y\r)\r)\fr{v^{p+1}}{p+1}\mrm dy.$$
On the left-hand side, we integrate by parts as before and then we apply Theorem \ref{ks}:
$$\fr{1}2\int_{\pa\Di}(\n v\cd\nu(y))^2\mrm d\s(y)\ge\fr{1}{4\pi}\l(\int_\Di\D v\mrm dy\r)^2=\fr{1}{4\pi}\l(\int_\Di|y|^{2\a}\wt K(y)v^p\mrm dy\r)^2\ge\fr{1}{4\pi C_0^2p^2}\l(p\int_\Di|y|^{2\a}\wt K(y)v^{p+1}\mrm dy\r)^2;$$
on the right-hand side, we argue as before exploiting the properties of $\wt K(y)$:
\beqy
&&2\int_\Di|y|^{2\a}\wt K(y)\fr{v^{p+1}}{p+1}\mrm dy+\int_\Di\l(2\a|y|^{2\a}\wt K(y)+|y|^{2\a}\l(\n\wt K(y)\cd y\r)\r)\fr{v^{p+1}}{p+1}\mrm dy\\
&\le&\ub{\l(2(1+\a)+\fr{\sup_{y\in\Di}\l|\n\wt K(y)\r|}{\inf_{y\in\Di}\wt K(y)}\r)}_{=:\fr{\rho_0}{8\pi C_0^2}}\int_\Di|y|^{2\a}\wt K(y)\fr{v^{p+1}}{p+1}\mrm dy\\
&\le&\fr{\rho_0}{4\pi C_0^2p^2}\l(p\int_\Di|y|^{2\a}\wt K(y)v^{p+1}\mrm dy\r).
\eeqy
Therefore we get
$$\fr{1}{4\pi C_0^2p^2}\l(p\int_\Di|y|^{2\a}\wt K(y)v^{p+1}\mrm dy\r)^2\le\fr{\rho_0}{4\pi C_0^2p^2}\l(p\int_\Di|y|^{2\a}\wt K(y)v^{p+1}\mrm dy\r),$$
namely
$$p\int_\O|x|^{2\a}K(x)u^{p+1}\mrm dx=p\int_\Di|y|^{2\a}\wt K(y)v^{p+1}\mrm dy\le\rho_0.$$
\epf\

\bpf[Proof of Theorem \ref{lious}]${}$\\
We apply the Riemann mapping to a solution $u=(u_1,\ds,u_N)$ to \eqref{liouseq} and we get a solution $v=(v_1,\ds,v_N)$ to $\l\{\bll-\D v_i=\sum_{j=1}^Na_{ij}\la_j|y|^{2\a_j}\wt K_j(y)e^{v_j}&\tx{in }\Di\\u_i=0&\tx{on }\pa\Di\earr\r.\q\q\q i=1,\ds,N$. Then, we apply to the latter Poho\v zaev identity, which gets
\beqy
&&\fr{1}2\sum_{i,j=1}^Na^{ij}\int_{\pa\Di}(\n v_i\cd\nu(y))(\n v_j\cd\nu(y))\mrm d\s(y)\\
&=&\sum_{i=1}^N\l(2\la_i\int_\Di|y|^{2\a_i}\wt K_i(y)e^{v_i}\mrm dy+\la_i\int_\Di\l(2\a_i|y|^{2\a_i}\wt K_i(y)+|y|^{2\a_i}\l(\n\wt K_i(y)\cd y\r)\r)e^{v_i}\mrm dy-\la\int_{\pa\Di}\wt K_i(y)e^{v_i}\mrm d\s(y)\r).
\eeqy
On the left-hand side we use the positivity of $A^{-1}$, as well as previous arguments, and we get
\beqy
&&\fr{1}2\sum_{i,j=1}^Na^{ij}\int_{\pa\Di}(\n v_i\cd\nu(y))(\n v_j\cd\nu(y))\mrm d\s(y)\\
&\ge&\fr{1}{4\pi}\sum_{i,j=1}^Na^{ij}\l(\int_{\pa\Di}(\n v_i\cd\nu(y))\mrm d\s(y)\r)\l(\int_{\pa\Di}(\n v_j\cd\nu(y))\mrm d\s(y)\r)\\
&=&\fr{1}{4\pi}\sum_{i,j=1}^Na^{ij}\l(\int_\Di\D v_i\mrm dy\r)\l(\int_\Di\D v_j\mrm dy\r)\\
&=&\fr{1}{4\pi}\sum_{i,j=1}^Na^{ij}\l(\la_i\int_\Di|y|^{2\a_i}\wt K_i(y)e^{v_i}\mrm dy\r)\l(\la_j\int_\Di|y|^{2\a_j}\wt K_j(y)e^{v_j}\mrm dy\r);
\eeqy
on the right-hand side we just use the same estimates as in Theorem \ref{liou} to each term and get
\beqy
&&\sum_{i=1}^N\l(2\la_i\int_\Di|y|^{2\a_i}\wt K_i(y)e^{v_i}\mrm dy+\la_i\int_\Di\l(2\a_i|y|^{2\a_i}\wt K_i(y)+|y|^{2\a_i}\l(\n\wt K_i(y)\cd y\r)\r)e^{v_i}\mrm dy-\la_i\int_{\pa\Di}\wt K_i(y)e^{v_i}\mrm d\s(y)\r)\\
&\le&C\sum_{i=1}^N\la_i\int_\Di|y|^{2\a_i}\wt K_i(y)e^{v_i}\mrm dy.
\eeqy
Finally, since $A^{-1}$ is positive definite, the relation we just found
$$\fr{1}{4\pi}\sum_{i,j=1}^Na^{ij}\l(\la_i\int_\Di|y|^{2\a_i}\wt K_i(y)e^{v_i}\mrm dy\r)\l(\la_j\int_\Di|y|^{2\a_j}\wt K_j(y)e^{v_j}\mrm dy\r)\le C\sum_{i=1}^N\la_i\int_\Di|y|^{2\a_i}\wt K_i(y)e^{v_i}\mrm dy$$
can only be satisfied if every integral belongs to a bounded region, namely
$$\la_i\int_\O|x|^{2\a_i}K_i(x)e^{u_i}\mrm dx=\la_i\int_\Di|y|^{2\a_i}\wt K_i(y)e^{v_i}\mrm dy\le\rho_0\q\q\q i=1,\ds,N$$

\epf\

\bpf[Proof of Theorem \ref{gen}]${}$\\
As before, we take a solution $u$ to \eqref{geneq} and apply the Riemann mapping theorem, thus getting a solution to $\l\{\bll-\D v=\wt W(y)F'(v)&\tx{in }\Di\\u=0&\tx{on }\pa\Di\earr\r.,$ with $\wt W(y)$ satisfying the same conditions as \eqref{cond}. Because of the properties of $\wt W(y)$ and $F(v)$, we have $v\in W^{2,q}(\Di)$ for some $q>1$, $|\D v|\le C|y|^{2\a}$ and $\l|\n\wt W(y)\r||y|\le C|y|^{2\a}$, therefore, as we argued in the proof of Theorem \ref{liou}, we can apply Theorem \ref{poho} to get:
$$\fr{1}2\int_{\pa\Di}(\n v\cd\nu(y))^2\mrm d\s(y)=2\int_\Di\wt W(y)F(v)\mrm dy+\int_\Di\l(\n\wt W(y)\cd y\r)F(v)\mrm dy-\int_{\pa\Di}\wt W(y)F(v)\mrm d\s(y).$$
Arguing as before we obtain
\beqy
&&\fr{1}{4\pi}\l(\int_\Di\wt W(y)F'(v)\mrm dy\r)^2\\
&\le&\fr{1}2\int_{\pa\Di}(\n v\cd\nu(y))^2\mrm d\s(y)\\
&=&2\int_\Di\wt W(y)F(v)\mrm dy+\int_\Di\l(\n_y\wt W(y)\cd y\r)F(v)\mrm dy-\int_{\pa\Di}\wt W(y)F(v)\mrm d\s(y)\\
&\le&(2+C)\int_\Di\wt W(y)F(v)\mrm dy\\
&\le&C(2+C)\int_\Di\wt W(y)F'(v)\mrm dy+C(2+C)\int_\Di\wt W(y)\mrm dy,
\eeqy
which means the mass must be uniformly bounded.
\epf\

\section*{Acknowledgments}\

The author wishes to thank Professor Daniele Bartolucci for the fruitful discussions concerning the topics of the paper.

\bibliography{nota2}

\begin{thebibliography}{10}

\bibitem{bp}
S.~Baraket and F.~Pacard.
\newblock Construction of singular limits for a semilinear elliptic equation in
  dimension {$2$}.
\newblock {\em Calc. Var. Partial Differential Equations}, 6(1):1--38, 1998.

\bibitem{bar}
D.~Bartolucci.
\newblock Existence and non existence results for supercritical systems of
  {L}iouville-type equations on simply connected domains.
\newblock {\em Calc. Var. Partial Differential Equations}, 53(1-2):317--348,
  2015.

\bibitem{bdm}
D.~Bartolucci, F.~De~Marchis, and A.~Malchiodi.
\newblock Supercritical conformal metrics on surfaces with conical
  singularities.
\newblock {\em Int. Math. Res. Not. IMRN}, (24):5625--5643, 2011.

\bibitem{barmal}
D.~Bartolucci and A.~Malchiodi.
\newblock An improved geometric inequality via vanishing moments, with
  applications to singular {L}iouville equations.
\newblock {\em Comm. Math. Phys.}, 322(2):415--452, 2013.

\bibitem{bat1}
L.~Battaglia.
\newblock Existence and multiplicity result for the singular {T}oda system.
\newblock {\em J. Math. Anal. Appl.}, 424(1):49--85, 2015.

\bibitem{bat0}
L.~Battaglia.
\newblock Variational aspects of singular {L}iouville systems.
\newblock {\em PhD thesis}, 2015.

\bibitem{bat2}
L.~Battaglia.
\newblock {$B_2$} and {$G_2$} {T}oda systems on compact surfaces: {A}
  variational approach.
\newblock {\em J. Math. Phys.}, 58(1):011506, 25, 2017.

\bibitem{bjmr}
L.~Battaglia, A.~Jevnikar, A.~Malchiodi, and D.~Ruiz.
\newblock A general existence result for the {T}oda {S}ystem on compact
  surfaces.
\newblock {\em Adv. Math.}, 285:937--979, 2015.

\bibitem{batmal}
L.~Battaglia and A.~Malchiodi.
\newblock Existence and non-existence results for the {$SU(3)$} singular {T}oda
  system on compact surfaces.
\newblock {\em J. Funct. Anal.}, 270(10):3750--3807, 2016.

\bibitem{car}
A.~Carlotto.
\newblock On the solvability of singular {L}iouville equations on compact
  surfaces of arbitrary genus.
\newblock {\em Trans. Amer. Math. Soc.}, 366(3):1237--1256, 2014.

\bibitem{cm}
A.~Carlotto and A.~Malchiodi.
\newblock Weighted barycentric sets and singular {L}iouville equations on
  compact surfaces.
\newblock {\em J. Funct. Anal.}, 262(2):409--450, 2012.

\bibitem{cl03}
C.-C. Chen and C.-S. Lin.
\newblock Topological degree for a mean field equation on {R}iemann surfaces.
\newblock {\em Comm. Pure Appl. Math.}, 56(12):1667--1727, 2003.

\bibitem{cl13}
C.-C. Chen and C.-S. Lin.
\newblock Mean field equation of {L}iouville type with singular data:
  topological degree.
\newblock {\em Comm. Pure Appl. Math.}, 68(6):887--947, 2015.

\bibitem{dln}
D.~G. de~Figueiredo, P.-L. Lions, and R.~D. Nussbaum.
\newblock A priori estimates and existence of positive solutions of semilinear
  elliptic equations.
\newblock {\em J. Math. Pures Appl. (9)}, 61(1):41--63, 1982.

\bibitem{dgip}
F.~De~Marchis, M.~Grossi, I.~Ianni, and F.~Pacella.
\newblock {$L^\infty$}-norm and energy quantization for the planar
  {L}ane-{E}mden problem with large exponent.
\newblock {\em preprint}, 2018.

\bibitem{dip}
F.~De~Marchis, I.~Ianni, and F.~Pacella.
\newblock Asymptotic profile of positive solutions of {L}ane-{E}mden problems
  in dimension two.
\newblock {\em J. Fixed Point Theory Appl.}, 19(1):889--916, 2017.

\bibitem{dl}
F.~De~Marchis and R.~L\'opez-Soriano.
\newblock Existence and non existence results for the singular {N}irenberg
  problem.
\newblock {\em Calc. Var. Partial Differential Equations}, 55(2):Art. 36, 35,
  2016.

\bibitem{dlr}
F.~De~Marchis, R.~L\'opez-Soriano, and D.~Ruiz.
\newblock Compactness, existence and multiplicity for the singular mean field
  problem with sign-changing potentials.
\newblock {\em J. Math. Pures Appl. (9)}, 115:237--267, 2018.

\bibitem{dkm}
M.~del Pino, M.~Kowalczyk, and M.~Musso.
\newblock Singular limits in {L}iouville-type equations.
\newblock {\em Calc. Var. Partial Differential Equations}, 24(1):47--81, 2005.

\bibitem{dgm}
S.~Deng, D.~Garrido, and M.~Musso.
\newblock Multiple blow-up solutions for an exponential nonlinearity with
  potential in {$\Bbb{R}^2$}.
\newblock {\em Nonlinear Anal.}, 119:419--442, 2015.

\bibitem{dm1}
S.~Deng and M.~Musso.
\newblock Bubbling solutions for an exponential nonlinearity in {$\Bbb{R}^2$}.
\newblock {\em J. Differential Equations}, 257(7):2259--2302, 2014.

\bibitem{dm2}
S.~Deng and M.~Musso.
\newblock Blow up solutions for a {L}iouville equation with {H}\'enon term.
\newblock {\em Nonlinear Anal.}, 129:320--342, 2015.

\bibitem{dja}
Z.~Djadli.
\newblock Existence result for the mean field problem on {R}iemann surfaces of
  all genuses.
\newblock {\em Commun. Contemp. Math.}, 10(2):205--220, 2008.

\bibitem{djamal}
Z.~Djadli and A.~Malchiodi.
\newblock Existence of conformal metrics with constant {$Q$}-curvature.
\newblock {\em Ann. of Math. (2)}, 168(3):813--858, 2008.

\bibitem{egp}
P.~Esposito, M.~Grossi, and A.~Pistoia.
\newblock On the existence of blowing-up solutions for a mean field equation.
\newblock {\em Ann. Inst. H. Poincar\'e Anal. Non Lin\'eaire}, 22(2):227--257,
  2005.

\bibitem{emp}
P.~Esposito, M.~Musso, and A.~Pistoia.
\newblock Concentrating solutions for a planar elliptic problem involving
  nonlinearities with large exponent.
\newblock {\em J. Differential Equations}, 227(1):29--68, 2006.

\bibitem{epw}
P.~Esposito, A.~Pistoia, and J.~Wei.
\newblock Concentrating solutions for the {H}\'enon equation in {$\Bbb R^2$}.
\newblock {\em J. Anal. Math.}, 100:249--280, 2006.

\bibitem{gnn}
B.~Gidas, W.~M. Ni, and L.~Nirenberg.
\newblock Symmetry and related properties via the maximum principle.
\newblock {\em Comm. Math. Phys.}, 68(3):209--243, 1979.

\bibitem{gp}
M.~Grossi and A.~Pistoia.
\newblock Multiple blow-up phenomena for the sinh-{P}oisson equation.
\newblock {\em Arch. Ration. Mech. Anal.}, 209(1):287--320, 2013.

\bibitem{gt}
M.~Grossi and F.~Takahashi.
\newblock Nonexistence of multi-bubble solutions to some elliptic equations on
  convex domains.
\newblock {\em J. Funct. Anal.}, 259(4):904--917, 2010.

\bibitem{ks}
N.~Kamburov and B.~Sirakov.
\newblock Uniform a priori estimates for positive solutions of the
  {L}ane-{E}mden equation in the plane.
\newblock {\em preprint}, 2018.

\bibitem{lr}
R.~L\'opez-Soriano and D.~Ruiz.
\newblock Prescribing the {G}aussian curvature in a subdomain of {$\Bbb{S}^2$}
  with {N}eumann boundary condition.
\newblock {\em J. Geom. Anal.}, 26(1):630--644, 2016.

\bibitem{mpw}
M.~Musso, A.~Pistoia, and J.~Wei.
\newblock New blow-up phenomena for {$SU(n+1)$} {T}oda system.
\newblock {\em J. Differential Equations}, 260(7):6232--6266, 2016.

\bibitem{pr}
A.~Pistoia and T.~Ricciardi.
\newblock Sign-changing tower of bubbles for a sinh-{P}oisson equation with
  asymmetric exponents.
\newblock {\em Discrete Contin. Dyn. Syst.}, 37(11):5651--5692, 2017.

\bibitem{rw1}
X.~Ren and J.~Wei.
\newblock On a two-dimensional elliptic problem with large exponent in
  nonlinearity.
\newblock {\em Trans. Amer. Math. Soc.}, 343(2):749--763, 1994.

\bibitem{rw2}
X.~Ren and J.~Wei.
\newblock Single-point condensation and least-energy solutions.
\newblock {\em Proc. Amer. Math. Soc.}, 124(1):111--120, 1996.

\end{thebibliography}
\bibliographystyle{abbrv}

\end{document}